\documentclass[11pt]{article}
\usepackage{amssymb}
\usepackage{amsmath}
\usepackage{amsfonts, latexsym}

\begin{document}

\title{A BMO theorem for $\epsilon$ distorted diffeomorphisms from  $\mathbb R^D$ to $\mathbb R^D$ with applications to  manifolds of speech and sound}
\author{ Steven B. Damelin \thanks{ Zentralblatt MATH, FIZ Karlsruhe – Leibniz Instituteemail: steve.damelin@gmail.com}\and Charles Fefferman, \thanks{Department of Mathematics; Fine Hall, Washington Road, Princeton NJ 08544-1000 USA,\, email: cf@math.princeton.edu}\and William Glover\thanks{Department of Music, Albany State University, Albany, GA.\, email: will25655@yahoo.com }}

\maketitle

\begin{abstract}
This paper deals with a  BMO Theorem for $\epsilon$ distorted diffeomorphisms from $\mathbb R^D$ to $\mathbb R^D$ with  applications to  manifolds of speech and sound. The work in this paper appears in the memoir \cite{SDam}.
\end{abstract}
\bigskip

\noindent 1991 AMS(MOS) Classification:  58C25, 42B35, 94A08, 94C30, 41A05, 68Q25, 30E05, 26E10, 68Q17.
\medskip


\section{Introduction}
\subsection{Music, Speech and Mathematics}
From the very beginning of time, mathematicians have been intrigued by the facinating connections which exist between music, speech and mathematics. Indeed, these connections were already in some subtle form in the writings of Gauss. The aim of this paper is to study estimates in measure for diffeomorphisms from $\mathbb R^D$ to $\mathbb R^D$, $D\geq 2$ of small distortion and provide an application to 
music and speech manifolds. 

\section{Preliminaries}
\setcounter{equation}{0}

Fix a dimension $D\geq 2$. We work in $\mathbb R^D$. We write $B(x,r)$ to denote the open ball in $\mathbb R^D$ with centre $x$ and radius $r$. We write $A$ to denote Euclidean motions on $\mathbb R^D$. A Euclidean motion may be orientation-preserving or orientation reversing.
We write $c$, $C$, $C'$ etc to denote constants depending on the dimension $D$. These expressions need not denote the same constant in different occurrences. For a $D\times D$ matrix, $M=(M_{ij})$, we write
$|M|$ to denote the Hilbert-Schmidt norm
\[
|M|=\left(\sum_{ij}|M_{ij}|^2\right)^{1/2}.
\]
Note that if $M$ is real and symmetric and if
\[
(1-\lambda)I\leq M\leq (1+\lambda)I
\]
as matrices, where $0<\lambda<1$, then
\begin{equation}
|M-I|\leq C\lambda.
\end{equation}
This follows from working in an orthonormal basis for which $M$ is diagonal. One way to understand the formulas above is to think of $\lambda$ as being close to zero. See also (2.6) below.

A function $f:\mathbb R^D\to \mathbb R$ is said to be BMO (Bounded mean oscillation )if there is a constant $K\geq 0$ such that, for every ball $B\subset \mathbb R^D$, there exists a real number $H_B$ such that
\begin{equation}
\frac{1}{{\rm vol}\, B}\int_{B}|f(x)-H_B|dx\leq K.
\end{equation}
The least such $K$ is denoted by $||f||_{{\rm BMO}}$.

In harmonic analysis, a function of bounded mean oscillation, also known as a BMO function, is a real-valued function whose mean oscillation is bounded (finite). The space of functions of bounded mean oscillation (BMO), is a function space that, in some precise sense, plays the same role in the theory of Hardy spaces, that the space of essentially bounded functions plays in the theory of $Lp$-spaces: it is also called a John-Nirenberg space, after Fritz John and Louis Nirenberg who introduced and studied it for the first time. See \cite{J,JN}.

The John-Nirenberg inequality asserts the following: Let $f\in BMO$ and let $B\subset \mathbb R^D$ be a ball. Then there exists a real number $H_B$ such that
\begin{equation}
{\rm vol}\left\{x\in B:\, |f(x)-H_B|>C\lambda ||f||_{BMO}\right\}\leq \exp(-\lambda){\rm vol}\,B,\, \lambda \geq 1.
\end{equation}

As a corollary of the John-Nirenberg inequality, we have
\begin{equation}
\left(\frac{1}{{\rm vol}\, B}\int_{B}|f(x)-H_B|^4dx\right)^{1/4}\leq C\lambda ||f||_{BMO}.
\end{equation}

There is nothing special about the 4th power in the above; it will be needed later.

The definition of BMO, the notion of the BMO norm, the John-Nirenburg inequality (2.3) and its corollary (2.4) carry through to the case of functions $f$ on $\mathbb R^D$ which take their values in the space
of $D\times D$ matrices. Indeed, we take $H_B$ in (2.2-2.4) to be a $D\times D$ matrix for such $f$. The matrix valued norms of (2.3-2.4) follow easily from the scalar case.

We will need some potential theory. If $f$ is a smooth function of compact support in $\mathbb R^D$, then we can write $\Delta^{-1}f$ to denote the convolution of $f$ with the Newtonian potential.
Thus, $\Delta^{-1}f$ is smooth and $\Delta(\Delta^{-1}f)=f$ on $\mathbb R^D$.

We will use the estimate:
\begin{equation}
\left\|\frac{\partial}{\partial x_i}\Delta^{-1}\frac{\partial}{\partial x_j}f\right\|_{L^2(\mathbb R^D)}
\leq C||f||_{L^2(\mathbb R^D)},\, ij=1,...,D
\end{equation}
valid for any smooth function $f$ with compact support. Estimate (2.5) follows by applying the Fourier transform.

We will work with a positive number $\varepsilon$. We always assume that $\varepsilon\leq {\rm min}(1,C)$.
An $\varepsilon$ distorted diffeomorphism of $\mathbb R^D$ is a one to one and onto diffeomorphism
$\Phi:\mathbb R^D\to \mathbb R^D$ such as
\[
(1-\varepsilon)I\leq (\Phi'(x))^{T}(\Phi'(x))\leq (1+\varepsilon)I
\]
as matrices. Thanks to (2.1), such $\Phi$ satisfy
\begin{equation}
\left|(\Phi'(x))^{T}(\Phi'(x))-I\right|\leq C\varepsilon.
\end{equation}

We end this section, with the following inequality from \cite{FD}:                                        \medskip

{\bf Approximation Lemma}\, Let $\Phi:\mathbb R^D\to \mathbb R^D$ be an $\varepsilon$ distorted diffeomorphism. Then, there exists an Euclidean motion $A$ such that
\begin{equation}
\left|\Phi(x)-A(x)\right|\leq C\varepsilon
\end{equation}
for all $x\in B(0, 10)$.

\section{An overdetermined system}
\setcounter{equation}{0}

We will need to study the following elemetary overdetermined system of partial differential equations.

\begin{equation}
\frac{\partial \Omega_i}{\partial x_j}+\frac{\partial \Omega_j}{\partial x_i}=f_{ij}, i,j=1,...,D
\end {equation}
on $\mathbb R^D$. Here, $\Omega_i$ and $f_{ij}$ are $C^{\infty}$ functions on $\mathbb R^D$.
A result concerning (3.1) we need is:
\medskip

{\bf PDE Theorem}\, Let $\Omega_1$,...,$\Omega_D$ and $f_{ij}$, $i,j=1,...,D$ be smooth functions on $\mathbb R^D$. Assume that (3.1) holds and suppose that
\begin{equation}
||f_{ij}||_{L^2(B(0,4))}\leq 1.
\end{equation}
Then, there exist real numbers $\Delta_{ij}$, $i,j=1,...,D$ such that
\begin{equation}
\Delta_{ij}+\Delta_{ji}=0,\, \forall i,j
\end{equation}
and
\begin{equation}
\left\|\frac{\partial \Omega_i}{\partial x_j}-\Delta_{ij}\right\|_{L^2(B(0,1))}\leq C.
\end{equation}
\medskip

{\bf Proof}\, From (3.1), we see at once that
\[
\frac{\partial \Omega_i}{\partial x_i}=\frac{1}{2}f_{ii}
\]
for each $i$. Now, by differentiating (3.1) with respect to $x_j$ and then summing on $j$, we see that
\[
\Delta \Omega_i +\frac{1}{2}\frac{\partial}{\partial x_i}\left(\sum_j f_{jj}\right)=\sum_j \frac{\partial f_{ij}}{\partial x_j}
\]
for each $i$.
Therefore, we may write
\[
\Delta \Omega_i=\sum_j \frac{\partial}{\partial x_j} g_{ij}
\]
for smooth functions $g_{ij}$ with
\[
||g_{ij}||_{L^2(B(0,4)}\leq C.
\]
This holds for each $i$. Let $\chi$ be a $C^{\infty}$ cutoff function on $\mathbb R^D$ equal to 1 on $B(0,2)$ vanishing outside $B(0,4)$ and satisfying $0\leq \chi\leq 1$ everywhere. Now let
\[
\Omega_i^{{\rm err}}=\Delta^{-1}\sum_j\frac{\partial}{\partial x_j}\left(\chi g_{ji}\right)
\]
and let
\[
\Omega_i^*=\Omega_i-\Omega_i^{err}.
\]
Then,
\begin{equation}
\Omega_i=\Omega_i^*+ \Omega_i^{err}
\end{equation}
each $i$.
\begin{equation}
\Omega_i^*
\end{equation}
is harmonic on $B(0,2)$ and
\begin{equation}
\left||\nabla \Omega_i^{{\rm err}}\right||_{L^2(B(0,2))}\leq C
\end{equation}
thanks to (2.5). By (3.1, 3.2, 3.5, 3.7), we can write
\begin{equation}
\frac{\partial \Omega_i^*}{\partial x_j}+\frac{\partial \Omega_j^*}{\partial x_i}=f_{ij}^*, i,j=1,...,D
\end{equation}
on $B(0,2)$ and with
\begin{equation}
\left||f_{ij}^*\right||_{L^2(B(0,2)}\leq C.
\end{equation}
From (3.6) and (3.8), we see that each $f_{ij}^*$ is a harmonic function on $B(0,2)$. Consequently,
(3.9) implies
\begin{equation}
sup_{B(0,1)}\left|\nabla f_{ij}^*\right|\leq C.
\end{equation}
From (3.8), we have for each $i,j,k$,
\begin{eqnarray}
&& \frac{\partial^2 \Omega_i^*}{\partial x_{j}\partial x_k}+ \frac{\partial^2 \Omega_k^*}{\partial x_{i}\partial x_j}=\frac{\partial f_{ik}^*}{\partial x_j}; \frac{\partial^2 \Omega_i^*}{\partial x_{j}\partial x_k}+ \frac{\partial^2 \Omega_j^*}{\partial x_{i}\partial x_k}=\frac{\partial f_{ij}^*}{\partial x_k} \\
&& \frac{\partial^2 \Omega_j^*}{\partial x_{i}\partial x_k}+ \frac{\partial^2 \Omega_k^*}{\partial x_{i}\partial x_j}=\frac{\partial f_{jk}^*}{\partial x_i}.
\end{eqnarray}
Now adding the first two equations above and subtracting the last, we obtain:
\begin{equation}
2\frac{\partial^2 \Omega_i^*}{\partial x_{j}\partial x_k}=\frac{\partial f_{ik}^*}{\partial x_j}+\frac{\partial f_{ij}^*}{\partial x_k}-\frac{\partial f_{jk}^*}{\partial x_i}
\end{equation}
on $B(0,1)$. Now from (3.10) and (3.13), we obtain the estimate
\begin{equation}
\left|\frac{\partial^2 \Omega_i^*}{\partial x_{j}\partial x_k}\right|\leq C
\end{equation}
on $B(0,1)$ for each $i,j,k$.
Now for each $i,j$, let
\begin{equation}
\Delta_{ij}^*=\frac{\partial \Omega_i^*}{\partial x_j}(0).
\end{equation}
By (3.14), we have
\begin{equation}
\left|\frac{\partial \Omega_i^*}{\partial x_j}-\Delta_{ij}^*\right|\leq C
\end{equation}
on $B(0,1)$ for each $i,j$. Recalling (3.5) and (3.7), we see that (3.16) implies that
\begin{equation}
\left\|\frac{\partial \Omega_i}{\partial x_j}-\Delta_{ij}^*\right\|_{L^2(B(0,1))}\leq C.
\end{equation}
Unfortunately, the $\Delta_{ij}^*$ need not satisfy (3.3). However, (3.1),  (3.2) and (3.17) imply the estimate
\[
\left|\Delta_{ij}^*+\Delta_{ji}^*\right|\leq C
\]
for each $i,j$. Hence, there exist real numbers $\Delta_{ij}$, $(i,j=1,...,D)$ such that
\begin{equation}
\Delta_{ij}+\Delta_{ji}=0
\end{equation}
and
\begin{equation}
\left|\Delta_{ij}^*-\Delta_{ij}\right|\leq C
\end{equation}
for each $i,j$. From (3.17) and (3.19), we see that
\begin{equation}
\left\|\frac{\partial \Omega_i}{\partial x_j}-\Delta_{ij}\right\|_{L^2(B(0,1))}\leq C
\end{equation}
for each $i$ and $j$.

Thus (3.18) and (3.20) are the desired conclusions of the Theorem. $\Box$

\section{A BMO Theorem}
\setcounter{equation}{0}

In this section, we prove the following:
\medskip

{\bf BMO Theorem 1}\, Let $\Phi:\mathbb R^D\to \mathbb R^D$ be an $\varepsilon$ diffeomorphism and let
$B\subset \mathbb R^D$ be a ball. Then, there exists $T\in O(D)$ such that
\begin{equation}
\frac{1}{{\rm vol}\, B}\int_{B}\left|\Phi'(x)-T\right|dx\leq C\varepsilon^{1/2}.
\end{equation}

{\bf Proof}\, Estimate (4.1) is preserved by translations and dilations. Hence we may assume that
\begin{equation}
B=B(0,1).
\end{equation}

Now we know that there exists an Euclidean motion
$A:\mathbb R^D\to \mathbb R^D$ such that
\begin{equation}
\left|\Phi(x)-A(x)\right|\leq C\varepsilon
\end{equation}
for $x\in B_{(0,10)}$.
Our desired conclusion (4.1) holds for $\Phi$ iff it holds for $A^{-1}o\Phi$ (with a different T). Hence, without loss of generality, we may assume that $A=I$. Thus, (4.3) becomes
\begin{equation}
\left|\Phi(x)-x\right|\leq C\varepsilon, x\in B(0,10).
\end{equation}
We set up some notation: We write the diffeomorphism $\Phi$ in coordinates by setting:
\begin{equation}
\Phi(x_1,...,x_D)=(y_1,...,y_D)
\end{equation}
where for each $i$, $1\leq i\leq D$,
\begin{equation}
y_i=\psi_i(x_1,...,x_D).
\end{equation}
First claim: For each $i=1,...,D,$
\begin{equation}
\int_{B(0,1)}\left|\frac{\partial \psi_i(x)}{\partial x_i}-1\right|\leq C\varepsilon.
\end{equation}

For this, for fixed $(x_2,...,x_D)\in B'$, we apply (4.4) to the points $x^{+}=(1,...,x_D)$ and $x^{-}=(1,...,x_D)$.
We have
\[
\left|\psi_1(x^+)-1\right|\leq C\varepsilon
\]
and
\[
\left|\psi_1(x^{-1})+1\right|\leq C\varepsilon.
\]
Consequently,
\begin{equation}
\int_{-1}^{1}\frac{\partial \psi_1}{\partial x_1}(x_1,...,x_D)dx_1\geq 2-C\varepsilon.
\end{equation}
On the other hand, since,
\[
\left(\psi'(x)\right)^{T}\left(\psi'(x)\right)\leq (1+\varepsilon)I,
\]
we have the inequality for each $i=1,...D$,
\[
\left(\frac{\partial \psi_i}{\partial x_i}\right)^2\leq 1+\varepsilon.
\]
Therefore,
\begin{equation}
\left|\frac{\partial \psi_i}{\partial x_i}\right|-1\leq \sqrt{1+\varepsilon}-1\leq \varepsilon.
\end{equation}
Set
\[
I^{+}=\left\{x_1\in [-1,1]:\, \frac{\partial \psi_1}{\partial x_1}(x_1,...,x_D)-1\leq 0\right\},
\]
\[
I^{-1}=\left\{x_1\in [-1,1]:\, \frac{\partial \psi_1}{\partial x_1}(x_1,...,x_D)-1\geq 0\right\},
\]
\[
\Delta^{+}=\int_{I^+}\left(\frac{\partial \psi_1}{\partial x_1}(x_1,...,x_D)-1\right)dx_1
\]
and
\[
\Delta^{-}=\int_{I^-}\left(\frac{\partial \psi_1}{\partial x_1}(x_1,...,x_D)-1\right)dx_1.
\]
The inequality (4.8) implies that $-\Delta^{-1}\leq C\varepsilon+\Delta^{+}$. The inequality (4.9) implies that
\[
\frac{\partial \psi_1}{\partial x_1}-1\leq C\varepsilon.
\]
Integrating the last inequality over $I^+$, we obtain $\Delta^+\leq C\varepsilon$. Consequently,
\begin{equation}
\int_{-1}^{1}\left|\frac{\partial \psi_1}{\partial x_1}(x_1,...,x_D)-1\right|dx_1=\Delta^+ -\Delta^-\leq C\varepsilon.
\end{equation}
Integrating this last equation over $(x_2,..., x_D)\in B'$ and noting that $B(0,1)\subset [-1,1]\times B'$, we conclude that
\[
\int_{B(0,1)}\left|\frac{\partial \psi_1}{\partial x_1}(x_1,...,x_D)-1\right|dx\leq C\varepsilon.
\]
Similarly, for each $i=1,...,D$, we obtain (4.7).
\medskip

Second claim: For each $i,j=1,...,D, i\neq j$,
\begin{equation}
\int_{B(0,1)}\left|\frac{\partial \psi_i(x)}{\partial x_j}\right|dx\leq C\sqrt{\varepsilon}.
\end{equation}

Since
\[
(1-\varepsilon)I\leq (\Phi'(x))^{T}(\Phi'(x))\leq (1+\varepsilon)I,
\]
we have
\begin{equation}
\sum_{i,j=1}^D \left(\frac{\partial \psi_i}{\partial x_j}\right)^2\leq (1+C\varepsilon)D.
\end{equation}

Therefore,
\[
\sum_{i\neq j}\left(\frac{\partial \psi_i}{\partial x_j}\right)^2
\leq C\varepsilon+\sum_{i=1}^D\left(1-\frac{\partial \psi_i}{\partial x_i}\right)
\left(1+\frac{\partial \psi_i}{\partial x_i}\right).
\]
Using (4.9) for $i$, we have $\left|\frac{\partial \psi_i}{\partial x_i}\right|+1\leq C$. Therefore,
\[
\sum_{i\neq j}\left(\frac{\partial \psi_i}{\partial x_j}\right)^2
\leq C\varepsilon+C\left\|\frac{\partial \psi_i}{\partial x_i}-1\right|.
\]

Now integrating the last inequality over the unit ball and using (4.7), we find that
\begin{equation} \int_{B(0,1)}\sum_{i\neq j}\left(\frac{\partial \psi_i}{\partial x_j}\right)^2dx\leq C\varepsilon+\int_{B(0,1)}
\left\|\frac{\partial \psi_i}{\partial x_i}-1\right|dx\leq C\varepsilon.
\end{equation}

Consequently, by the Cauchy-Schwartz inequality, we have
\[
\int_{B(0,1)}\sum_{i\neq j}\left|\frac{\partial \psi_i}{\partial x_j}\right|dx\leq C\sqrt{\varepsilon}.
\]

Third claim:

\begin{equation}
\int_{B(0,1)}\left|\frac{\partial \psi_i}{\partial x_i}\right|dx\leq C\sqrt{\varepsilon}.
\end{equation}

Since,
\[
\int_{B(0,1)} \left(\frac{\partial \psi_i}{\partial x_i}-1\right)^2dx\leq \int_{B(0,1)}
\left|\frac{\partial \psi_i}{\partial x_i}-1\right|\left|\frac{\partial \psi_i}{\partial x_i}+1\right|dx,
\]
using (4.7) and $\left|\frac{\partial \psi_i}{\partial x_i}\right|\leq 1+C\varepsilon$, we obtain
\[
\int_{B(0,1)}\left(\frac{\partial \psi_i}{\partial x_i}\right)^2dx\leq C\varepsilon.
\]
Thus, an application of Cauchy Schwartz, yields (4.15).
\medskip

Final claim: By the Hilbert Schmidt definition, we have
\begin{eqnarray*}
&& \int_{B(0,1)}|\Psi'(x)-I|dx=\int_{B(0,1)}\left(\sum_{i,j=1}^D\left(\frac{\partial \psi_i}{\partial x_j}-
\delta_{ij}\right)^2\right)^{1/2} \\
&& \leq \int_{B(0,1)}\sum_{i,j=1}^D\left|\frac{\partial \psi_i}{\partial x_j}-\delta_{ij}\right|dx.
\end{eqnarray*}
The estimate (4.11) combined with (4.15) yields:
\[
\int_{B(0,1)}\left|\Phi'(x)-I\right|dx \leq C\varepsilon^{1/2}.
\]
Thus we have proved (4.1) with $T=I$. The proof of the BMO Theorem 1 is complete. $\Box$
\medskip

{\bf Corollary}\, Let $\Phi:\mathbb R^D\to \mathbb R^D$ be an $\varepsilon$-distorted diffeomorphism. For each, ball $B\subset \mathbb R^D$, there exists $T_B\in O(D)$, such that
\[
\left(\frac{1}{{\rm vol}\, B}\int_{B}\left|\Phi'(x)-T\right|^4 dx\right)^{1/4}\leq C\varepsilon^{1/2}.
\]
The proof follows from the first BMO Theorem just proved and the John Nirenberg inequality. (See (2.4).
$\Box$.

\section{A Refined BMO Theorem}
\setcounter{equation}{0}

We prove:
\medskip

{\bf BMO Theorem 2}\, Let $\Phi:\mathbb R^D\to \mathbb R^D$ be an $\varepsilon$ diffeomorphism and let
$B\in \mathbb R^D$ be a ball. Then, there exists $T\in O(D)$ such that
\begin{equation}
\frac{1}{{\rm vol}\, B}\int_{B}\left|\Phi'(x)-T\right|dx\leq C\varepsilon.
\end{equation}

{\bf Proof}\,  We may assume without loss of generality that
\begin{equation}
B=B(0,1).
\end{equation}
We know that there exists $T_B^*\in O(D)$ such that
\[
\left(\int_{B}|\Phi'(x)-T^*_{B}|^4 dx\right)^{1/4}\leq C\varepsilon^{1/2}.
\]
Our desired conclusion holds for $\Phi$ iff it holds for $(T_{B}^{*})^{-1}o\Phi$. Hence without loss of generality, we may assume that $T_B^{*}=I$. Thus we have
\begin{equation}
\left(\int_{B}|\Phi'(x)-I)|^4 dx\right)^{1/4}\leq C\varepsilon^{1/2}.
\end{equation}
Let
\begin{equation}
\Omega(x)=\left(\Omega_1(x),\Omega_2(x),....,\Omega_{D}(x)\right)=\Phi(x)-x,\, x\in \mathbb R^{D}.
\end{equation}
Thus (5.3) asserts that
\begin{equation}
\left(\int_{B(0,1)}\left|\nabla\Omega(x)\right|^4dx\right)^{1/4}\leq C\varepsilon^{1/2}.
\end{equation}
We know that
\begin{equation}
\left|(\Phi'(x))^{T}\Phi'(x)-I\right|\leq C\varepsilon,\, x\in \mathbb R^{D}.
\end{equation}
In coordinates, $\Phi'(x)$  is the matrix
$\left(\delta_{ij}+\frac{\partial \Omega_i(x)}{\partial x_j}\right)$, hence
$\Phi'(x)^{T}\Phi'(x)$  is the matrix whose $ij$th entry is
\[
\delta_{ij}+\frac{\partial \Omega_j(x)}{\partial x_i}+\frac{\partial \Omega_i(x)}{\partial x_j}
+\sum_{l}\frac{\partial \Omega_l(x)}{\partial x_i}\frac{\partial \Omega_l(x)}{\partial x_j}.
\]
Thus (5.6) says that
\begin{equation}
\left|\frac{\partial \Omega_j}{\partial x_i}+\frac{\partial \Omega_i}{\partial x_j}
+\sum_{l}\frac{\partial \Omega_l}{\partial x_i}\frac{\partial \Omega_l}{\partial x_j}\right|\leq C
\varepsilon
\end{equation}
on $\mathbb R^{D},\, i,j=1,...,D.$
Thus, we have from (5.5), (5.7) and the Cauchy Schwartz inequality the estimate
\[
\left\|\frac{\partial \Omega_i}{\partial x_j}+\frac{\partial \Omega_j}{\partial x_i}\right\|_{L^2(B(0,10))}\leq C\varepsilon.
\]
By the PDE Theorem, there exists, for each $i,j$, an antisymmetric matrix $S=(S)_{ij}$, such that
\begin{equation}
\left\|\frac{\partial \Omega_i}{\partial x_j}-S\right\|_{L^2(B(0,1))}\leq C\varepsilon.
\end{equation}
Recalling (5.4), this is equivalent to
\begin{equation}
\left\|\Phi'-(I+S)\right\|_{L^2(B(0,1))}\leq C\varepsilon.
\end{equation}
Note that (5.5) and (5.8) show that
\[
|S|\leq C\varepsilon^{1/2}
\]
and thus,
\[
\left|\exp(S)-(I+S)\right|\leq C\varepsilon.
\]
Hence, (5.9) implies via Cauchy Schwartz.
\begin{equation}
\int_{B(0,1)}\left|\Phi'(x)-\exp(S)(x)\right|dx \leq C\varepsilon^{1/2}.
\end{equation}
This implies the result because $S$ is antisymmetric which means that $\exp(S)\in O(D)$. $\Box$.

\section{A BMO Theorem for Diffeomorphisms of Small Distortion}
\setcounter{equation}{0}

In this section, we  prove the following theorem.
\medskip

{\bf Theorem}\, Let $\Phi:\mathbb R^D\to \mathbb R^D$ be an $\varepsilon$ distorted diffeomorphism. Let $B\subset \mathbb R^D$ be a ball. Then, there exists $T_B\in O(D)$ such that for every $\lambda\geq 1$,
\begin{equation}
{\rm vol}\left\{x\in B:\, |\Phi'(x)-T_B|>C\lambda \varepsilon\right\}\leq \exp(-\lambda){\rm vol}\,(B).
\end{equation}

Moreover, the result (6.1) is sharp in the sense of small volume if one takes a slow twist defined as follows:
For $x\in \mathbb R^D$, let $S_x$ be the block-diagonal matrix
\[
\left(
\begin{array}{llllll}
D_1(x) & 0 & 0 & 0 & 0 & 0 \\
0 & D_2(x) & 0 & 0 & 0 & 0 \\
0 & 0 & . & 0 & 0 & 0 \\
0 & 0 & 0 & . & 0 & 0 \\
0 & 0 & 0 & 0 & . & 0 \\
0 & 0 & 0 & 0 & 0 & D_r(x)
\end{array}
\right)
\]
where, for each $i$, either $D_i(x)$ is the $1\times 1$ identity matrix or else
\[
D_i(x)=\left(
\begin{array}{ll}
\cos f_i(|x|) & \sin f_i(|x|) \\
-\sin f_i(|x|) & \cos f_i(|x|)
\end{array}
\right)
\]
for a function $f_i$ of one variable.

Now define for each $x\in \mathbb R^D$,
$\Phi(x)=\Theta^{T}S_x(\Theta x)$ where $\Theta$
is any fixed matrix in $SO(D)$.
One checks that $\Phi$ is $\varepsilon$-distorted, provided for each $i$,
$t|f_i'(t)|<c\varepsilon$ for all $t\in [0,\infty)$.
\medskip

{\bf Proof}\, The theorem follows from the BMO Theorem 2 and the John Nirenburg inequality. The sharpness can be easily checked.
$\Box$.

\section{On the Approximate and Exact Allignment of Data in Euclidean Space, Speech and Music Manifolds}

\subsection{Approximate and Exact Allignment of Data}

The following is a classical question in Euclidean Geometry, see for example \cite{WW}:
Suppose we are  that given two sets of distinct data points in Euclidean $D\geq 2$ space, say from two manifolds.
We do not know what the manifolds are apriori but we do know that the pairwise distances between the points are equal. Does there exist an isometry (distance preserving map) that maps the one set of points 1-1 onto the other. This is a fundamental question in data analysis for most often, we are only given sampled function points from two usually unknown manifolds and we seek to know what can be said about the manifolds themselves. A typical example might be a face recognition problem where all we have is multiple finite images of people's faces from various views.

An added complication in the above question is that in general, we are not given exact distances between function value points. We have noise and so we need to demand that instead of the
pairwise distances being equal, they should be " close" in some reasonable metric. It is well known, see \cite{WW} that any isometry of a subset of a Euclidean space into the space can be extended to an isometry of that Euclidean space onto itself. Some results on almost isometries in Euclidean spaces can be found for example in \cite{J}, \cite{ATV}.
Some results on almost isometries in Euclidean spaces can be found for example in \cite{J} and \cite{ATV}.
\medskip

Consider  the following theorems (\cite{SDam}) which tell us alot about how to handle manifold identification when the point set function values given are not exactly equal but are close.
\medskip

{\bf Theorem}\, Given $\varepsilon>0$ and $k\geq 1$, there exists $\delta>0$ such that the following holds. Let $y_1,...y_k$ and $z_1,...,z_k$ be points in $\mathbb R^D$. Suppose
\[
(1+\delta)^{-1}\leq \frac{|z_i-z_j|}{|y_i-y_j|}\leq 1+\delta,\, i\neq j.
\]
Then, there exists a Euclidean motion $\Phi_0:x\to Tx+x_0$ such that
\[
|z_i-\Phi_0(y_i)|\leq \varepsilon {\rm diam}\left\{y_1,...,y_k\right\}
\]
for each $i$. If $k\leq D$, then we can take $\Phi_0$ to be a proper Euclidean motion on $\mathbb R^{D}$.
\medskip

{\bf Theorem}\, Let $\varepsilon>0$, $D\geq 1$ and $1\leq k\leq D$.
Then there exists $\delta>0$ such that the following holds: Let $E:=y_1,...y_k$ and $E':=z_1,...z_k$ be distinct points in $\mathbb R^D$. Suppose that
\[
(1+\delta)^{-1}\leq \frac{|z_i-z_j|}{|y_i-y_j|}\leq (1+\delta),\,  1\leq i,j\leq k,\, i\neq j.
\]
Then there exists a diffeomorphism, 1-1 and onto map $\Psi:\mathbb R^D\to \mathbb R^D$ with
\[
(1+\varepsilon)^{-1}\leq \frac{|\Psi(x)-\Psi(y)|}{|x-y|}\leq (1+\varepsilon),\, x,y\in \mathbb R^D,\, x\neq y
\]
satisfying
\[
\Psi(y_i)=z_i,\, 1\leq i\leq k.
\]

Given the two theorems above, we now need to ask ourselves. Can we take, in any particular data application, a smooth $\varepsilon$ distortion and 
approximate it by an element of $O(D)$. Clearly this is very important. We understand that the results of this paper  tell us that at least in measure,
the derivative of a smooth $\varepsilon$ distortion may be well approximated by an element of $O(D)$. 
\subsection{Speech and Music manifolds}

Recently, see for example, \cite{DM} and the references cited therein there has been much interest in geometrically motivated dimensionality reduction algorithms. The reason for this is that these algorithms exploit low dimensional manifold structure in certain natural datasets to reduce dimensionality while preserving categorical content. In \cite{JN1}, the authors motivated the existence of low dimensional music and speech manifold structure to the existence of  certain rigid motions approximating smooth distortions of voice and speech sounds maps. The theorems proved in this paper and in \cite{FD} provide a fascinating insight into these very interesting questions.

 \end{document}